\documentclass{article}

\usepackage{amsmath}
\usepackage{amssymb}
\usepackage{theorem}
\usepackage{xypic}
\xyoption{all}

\newtheorem{thm}{Theorem}[section]
\newtheorem{cor}[thm]{Corollary}
\newtheorem{lem}[thm]{Lemma}
\newtheorem{prop}[thm]{Proposition}

\theorembodyfont{\upshape}
\newtheorem{defn}[thm]{Definition}

\newtheorem{rmk}[thm]{Remark}

\newtheorem{ex}[thm]{Example}

\newenvironment{pf}{\begin{trivlist}\item[]{\sc Proof.}}%
            {\nolinebreak $\Box$ \end{trivlist}}

\newenvironment{proof}{\begin{trivlist}\item[]{\sc Proof.}}%
            {\nolinebreak $\Box$ \end{trivlist}}

\newcommand{\toto}{\rightrightarrows}
\renewcommand{\c}{{\mathfrak c}}
\newcommand{\mult}{\mathop{\rm mult}\nolimits}
\newcommand{\vir}{{\rm vir}}

\newcommand{\ol}{\overline}

\newcommand{\A}{\mathbb A}

\newcommand{\cc}{\mathbb C}
\newcommand{\qq}{\mathbb Q}

\renewcommand{\L}{\mathfrak L}
\newcommand{\zz}{\mathbb Z}

\newcommand{\del}{\partial}

\DeclareMathOperator\id{id}

\DeclareMathOperator\spec{Spec}

\renewcommand{\O}{\mathcal O}

\newcommand{\EE}{\mathfrak E}
\newcommand{\CC}{\mathfrak C}

\newcommand{\cok}{\mathop{\rm cok}\nolimits}

\newcommand{\rk}{\mathop{\rm rk}\nolimits}

\newcommand{\Con}{\mathop{\rm Con}\nolimits}
\newcommand{\Ch}{\mathop{\rm Ch}\nolimits}
\newcommand{\Eu}{\mathop{\rm Eu}\nolimits}

\newcommand{\sheafhom}{\mathop{\mathcal{H}om}\nolimits}

\newcommand{\dotarrow}{{\SelectTips{cm}{}\xymatrix@1@C=10pt{{}\ar@{.>}[r]&{}}}} 

\newcommand{\noprint}[1]{}

\def\Label#1{\label{#1}{\tt [#1]}\phantom{h}}

\def\Label{\label}

\author{K. Behrend}
\title{Donaldson-Thomas Type Invariants via Microlocal
  Geometry}
\date{December 1, 2005}

\begin{document}
\sloppy

\maketitle

\abstract{We prove that Donaldson-Thomas type invariants are equal to
weighted Euler characteristics of their moduli spaces. In particular,
such invariants depend only on the scheme structure of the moduli
space, not the symmetric obstruction theory used to define them.  We
also introduce new invariants generalizing Donaldson-Thomas type
invariants to moduli problems with open moduli space. These are useful
for computing Donaldson-Thomas type invariants over stratifications. }

\eject
\tableofcontents
\eject

\section*{Introduction}
\addcontentsline{toc}{section}{Introduction}

\subsubsection{Donaldson-Thomas type invariants}

Donaldson-Thomas invariants \cite{RT} are the virtual counts of stable
sheaves (with fixed determinant) on Calabi-Yau threefolds.
Heuristically, the Donaldson-Thomas moduli space is the critical set
of the holomorphic Chern-Simons functional and the Donaldson-Thomas
invariant is a holomorphic analogue of the Casson invariant.

Recently \cite{MNOP}, Donaldson-Thomas invariants for sheaves of rank
one have been conjectured to have deep connections with Gromov-Witten
theory of Calabi-Yau threefolds.  They are supposed to encode the
integrality properties of such Gromov-Witten invariants, for example.

Mathematically, Donaldson-Thomas invariants are constructed as follows
(see \cite{RT}).  Deformation theory gives rise to a {\em perfect
obstruction theory }\cite{BF} (or a {\em tangent-obstruction complex
}in the language of \cite{LT}) on the moduli space of stable sheaves
$X$.
As Thomas points out in \cite{RT}, the {\em obstruction sheaf }is
equal to $\Omega_X$, the sheaf of K\"ahler differentials, and hence
the tangents $T_X$ are {\em dual }to the obstructions.  This expresses
a certain symmetry of the obstruction theory and is a mathematical
reflection of the heuristic that views $X$ as the critical locus of
a holomorphic functional.

Associated to the perfect obstruction theory is the {\em virtual
fundamental class}, an element of the Chow group $A_\ast(X)$ of
algebraic cycles modulo rational equivalence on $X$.  One implication
of the symmetry of the obstruction theory is the fact that the virtual
fundamental class $[X]^\vir$ is of degree zero. It can hence be
integrated over the proper space of stable sheaves to a number, the
{\em Donaldson-Thomas invariant }or {\em virtual count }of $X$
$$\#^\vir(X)=\int_{[X]^\vir}1\,.$$

We take the point of view that the symmetry of the obstruction theory
is the distinguishing feature of Donaldson-Thomas invariants, and call
any virtual count of a proper scheme with symmetric obstruction theory
a {\em Donaldson-Thomas type }invariant.

\subsubsection{Euler characteristics and $\nu_X$}

Note that if the moduli space $X$ is smooth, the obstruction sheaf
$\Omega_X$ is a bundle, so the virtual fundamental class is
the top Chern class $e(\Omega_X)$ and so the virtual count
is, up to a sign, the Euler characteristic of $X$:
\begin{equation*}
\#^\vir(X)=\int_{[X]}e(\Omega_X)=(-1)^{\dim X}\chi(X)\,,
\end{equation*}
We will generalize this formula to arbitrary (embeddable) schemes.

More precisely, we will construct on any scheme $X$ over $\cc$ in a
canonical way a constructible function $\nu_X:X\to\zz$ (depending only
on the scheme structure of $X$), such that if $X$ is proper and
embeddable we have 
\begin{equation}\Label{euler-main}
\#^\vir(X)=\chi(X,\nu_X)=\sum_{n\in\zz}n\,\chi\{\nu_X=n\}\,,
\end{equation}
for {\em any }symmetric obstruction theory  on $X$ with associated
Donaldson-Thomas type invariant $\#^\vir(X)$. 

As consequences of this result we obtain:

\smallskip\noindent$\bullet$
Donaldson-Thomas type invariants depend only on the scheme
  structure of the underlying moduli space, not on the symmetric
  obstruction theory used to define them.

\smallskip\noindent$\bullet$
 Even if $X$ is not proper, and so the virtual count does not
  make sense as an integral, we
  can consider the weighted Euler characteristic 
$$\tilde\chi(X)=\chi(X,\nu_X)$$
as a
  substitute for the virtual count.  This generalizes Donaldson-Thomas
  type invariants to the case of non-proper moduli space $X$. It also
  makes Donaldson-Thomas 
invariants accessible to arguments involving stratifying the moduli
space $X$. For an application, see~\cite{bf5}.

\smallskip\noindent$\bullet$
The value of $\nu_X$ at the point $P\in X$ should be considered as
the contribution of the point $P$ to the virtual count of $X$.

\smallskip
Some of the fundamental properties of $\nu_X$ are:

\smallskip\noindent$\bullet$ At smooth points $P$ of $X$ we have
$\nu_X(P)=(-1)^{\dim X}$.

\smallskip\noindent$\bullet$
 If $f:X\to Y$ is \'etale, then $f^\ast\nu_Y=\nu_X$.  Thus
  $\nu_X(P)$ is an invariant of the singularity of $X$ at the point
  $P$. 

\smallskip\noindent$\bullet$ Multiplicativity: $\nu_{X\times
  Y}(P,Q)=\nu_X(P)\nu_Y(Q)$. 

\smallskip\noindent$\bullet$ If $X$ is the critical scheme of a
regular function $f$ on a smooth
scheme $M$, i.e. $X=Z(df)$, then 
$$\nu_X(P)=(-1)^{\dim M}\big(1-\chi(F_P)\big)\,,$$
where $F_P$ is the Milnor fibre, i.e., the intersection of a nearby
fibre of $f$ with a small ball in $M$ centred at $P$.

\smallskip
Thus, if $X$ is the Donaldson-Thomas moduli space of stable sheaves,
one can, heuristically, think of $\nu_X$ as the Euler characteristic
of the perverse sheaf of vanishing cycles of the holomorphic
Chern-Simons functional.

Note that the existence of a symmetric obstruction theory on $X$ puts
strong restrictions on the singularities $X$ may have.  For example,
reduced local complete intersection singularities are excluded.  Thus
it is not clear how useful or significant $\nu_X$ is on general
schemes, which do not admit symmetric obstruction theories.

\subsubsection{Microlocal geometry}

Embed $X$ into a smooth scheme $M$.  Then we have a commutative diagram
\begin{equation}\Label{pretty}
\vcenter{\xymatrix{
Z_\ast(X)\rto^{\Eu}_{\sim}\drto_{c_0^M} &
\Con(X)\dto^{c_0^{SM}}\rto_{\sim}^{\Ch} & 
  \L_X(\Omega_M) \dlto^{0^!}\\
& A_0(X) & }}
\end{equation}
where the two horizontal arrows are isomorphisms. Here
$\Eu:Z_\ast(X)\to\Con(X)$ is MacPherson's local Euler obstruction
\cite{csm}, which maps algebraic cycles to $\zz$-valued constructible
functions and $\Ch:\Con(X)\to\L_X(\Omega_M)$ maps a constructible function to
its characteristic cycle, which is a conic Lagrangian cycle on
$\Omega_M$ supported inside $X$. The maps to $A_0(X)$ are the degree
zero  Chern-Mather class, the degree zero 
Schwartz-MacPherson Chern class, and the intersection with the zero
section, respectively. (Of course, the left part of the diagram exists
without the  embedding into $M$.) 

Now, given a symmetric obstruction theory on $X$, the cone of curvilinear
obstructions $cv\hookrightarrow ob=\Omega_X$, pulls back to a cone in 
$\Omega_M|_X$ via the epimorphism $\Omega_M|X\to \Omega_X$.  Via the
embedding $\Omega_M|_X\hookrightarrow \Omega_M$ we obtain a conic
subscheme $C\hookrightarrow \Omega_M$, the {\em obstruction cone }for
the embedding $X\hookrightarrow M$.  The virtual fundamental class is
$[X]^\vir=0^![C]$. 

The key fact is that $C$ is {\em Lagrangian}. Because of this, there
exists a unique constructible function $\nu_X$ on $X$ such that
$\Ch(\nu_X)=[C]$ and $c_0^{SM}(\nu_X)=[X]^\vir$.  Then
(\ref{euler-main}) follows by a simple application of MacPherson's
theorem~\cite{csm} (or equivalently from the microlocal index theorem of
Kashiwara~\cite{kash-schap}). 

The cycle $\c_X$ such that $\Eu(\c_X)=\nu_X$ is easily written
down. It can be thought of as the (signed) support of the intrinsic
normal cone of $X$.

The class $\alpha_X=c^M(\c_X)=c^{SM}(\nu_X)$, whose degree zero
component is the virtual fundamental class of any symmetric
obstruction theory on $X$, was  introduced by
Aluffi~\cite{Alu} (although with a different sign)  and
we call it therefor the {\em Aluffi class }of $X$.

We do not know if every scheme admitting a symmetric obstruction
theory can locally be written as the critical locus of a regular
function on a smooth scheme. This limits the usefulness of the above
formula for $\nu_X(P)$ in terms of the Milnor fibre.  Hence we provide
an alternative formula~(\ref{spirit}), similar in spirit, which always
applies.

If $\mathcal{M}$ is a regular holonomic $\mathcal{D}$-module on $M$
whose characteristic cycle is $[C]$, then
$$\nu_X(P)=\sum_i (-1)^i \dim_\cc H^i_{\{P\}}(X,\mathcal{M}_{DR})\,,$$
for any point $P\in M$.  Here $H^i_{\{P\}}$ denotes cohomology with
supports in the subscheme $\{P\}\hookrightarrow M$ and
$\mathcal{M}_{DR}$ denotes the perverse sheaf associated to
$\mathcal{M}$ via the Riemann-Hilbert correspondence, as incarnated,
for example, by the De~Rham complex $\mathcal{M}_{DR}$.  It would be
interesting to construct $\mathcal{M}$ or $\mathcal{M}_{DR}$ in
special cases, for example the moduli space of sheaves.  Maybe, as
opposed to $[C]$ and $[X]^\vir$, this more subtle data could actually
depend on the symmetric obstruction theory.

\subsubsection{Conventions}

We will always work over the field of complex numbers $\cc$.  All
schemes and algebraic stacks we consider are of finite type (over
$\cc$).  The  relevant facts about algebraic cycles and intersection
theory on stacks can be found in~\cite{vistoli} and~\cite{Kresch}.

We will often have to assume that our Deligne-Mumford stacks have the
resolution property or are embeddable into smooth stacks.  We
therefore consider quasi-projective Deligne-Mumford stacks
(see~\cite{quasi}): 

\begin{defn}[Kresch]
A separated Deligne-Mumford stack $X$, of finite type over $\cc$, with
quasi-projective coarse moduli space  is
called {\bf quasi-projective}, if any of the following equivalent
conditions are satisfied:

(i) $X$ has the resolution property, i.e., every coherent
$\O_X$-module is a quotient of a locally free coherent $\O_X$-module,

(ii) $X$ admits a finite flat cover $Y\to X$, where
$Y$ is a quasi-projective scheme,

(iii) $X$ is isomorphic to a quotient stack $[Y/GL_n]$, for some $n$,
where $Y$ is a 
quasi-projective scheme with a linear $GL_n$-action,

(iv) $X$ can be embedded as a locally closed substack into a smooth
separated Deligne-Mumford stack of finite type with projective coarse
moduli space.
\end{defn}

For  $\zz$-valued functions $f$, $g$ on sets $X$, $Y$, respectively,
we denote by $f\boxdot g$ the function on $X\times Y$ defined by
$(f\boxdot g)(x,y)=f(x)g(y)$, for all $(x,y)\in X\times Y$.

We will often use homological notation for complexes.
This means that $E_n=E^{-n}$, for a complex $\ldots\to
E^{i}\to E^{i+1}\to\ldots$ in some abelian category.

For a complex of sheaves $E$, we denote the cohomology sheaves by
$h^i(E)$. 

Let us recall a few sign conventions: If
$E=[E_1\stackrel{\alpha}{\longrightarrow} E_0]$ is a complex 
concentrated in the interval $[-1,0]$, then the dual complex
$E^\vee=[E_0^\vee\stackrel{-\alpha^\vee}{\longrightarrow} E_1^\vee]$
is a complex concentrated in the interval $[0,1]$. Thus the shifted
dual $E^\vee[1]$ is given by
$E^\vee[1]=[E_0^\vee\stackrel{\alpha^\vee}{\longrightarrow} E_1^\vee]$
and concentrated, again, in the interval $[-1,0]$. 

If $\theta:E\to F$ is a homomorphism of complexes concentrated in the
interval $[-1,0]$, such that
$\theta=(\theta_1,\theta_0)$, then the shifted
dual $\theta^\vee[1]:F^\vee[1]\to E^\vee[1]$ is given by
$\theta^\vee[1]=(\theta_0^\vee,\theta_1^\vee)$.

\subsubsection{Acknowledgements}

I would like to thank the following for helpful discussions: Paolo
Aluffi, Paul Bressler, Jim Bryan, Barbara Fantechi, Minhyong Kim and
Bob MacPherson. Special thanks are due to Yong-Geun Oh, who suggested
a microlocal approach to the problem and to Richard Thomas, who
suggested the key fact, that the obstruction cone is Lagrangian.

I would like to thank Sheldon Katz and the University of Illinois at
Urbana-Champaign for organizing the workshop on Donaldson-Thomas
invariants in March 2005, where I was able to speak on related work
and comments by participants were instrumental in bringing about this
paper.

Finally, I would like to thank Bumsig Kim and the Korea Institute for
Advanced Study for the hospitality during my stay in Seoul, where this
work was completed.  

This research was supported by a Discovery Grant from the Natural Sciences
and Engineering Research Council of Canada.

\eject
\section{A few invariants of schemes and stacks}

\subsection{The signed support of the intrinsic normal cone $\c_X$}

Let $X$ be a scheme.  Suppose $X$ is embedded as a closed subscheme of
a smooth scheme $M$. Consider the normal cone $C=C_{X/M}$ and its
projection $\pi:C\to X$. Define 
the cycle $\c_{X/M}$ on $X$ by
$$\c_{X/M}=\sum_{C'} (-1)^{\dim\pi(C')}\mult(C')\pi(C')\,.$$
The sum is over all irreducible components $C'$ of $C$.  By $\pi(C')$
we denote the irreducible closed subset (prime cycle) of $X$ obtained
as the image 
of $C'$ under $\pi$.  Alternatively, we can define $\pi(C')$ as the
(set-theoretic) 
intersection of $C'$ with the zero section of $C\to X$. The
multiplicity of the component $C'$ in the fundamental cycle $[C]$ of
$C$ is denoted by $\mult(C')$.  Hence $\mult(C')$ is the length of $C$ at
the generic point of $C'$. Note that even though
$[C]=\sum_{C'}\mult(C')C'$ is an effective cycle of homogeneous degree
$\dim M$, the cycle $\c_{X/M}$ is neither effective nor homogeneous. 

\begin{prop}
Let $X$ be a Deligne-Mumford stack. 
There is a unique (integral) cycle $\c_X$ on $X$ with the property that for any
\'etale map $U\to X$, and any closed embedding $U\to M$ of $U$ into a
smooth scheme $M$ we have
$$\c_X|_U=\c_{U/M}\,.$$
\end{prop}
\begin{pf}
Suppose given a commutative diagram of schemes 
$$\xymatrix{Y\rto\dto & N\dto\\X\rto & M}$$
where $Y\to N$ and $X\to M$ are closed embeddings, $Y\to X$ is \'etale
and $M\to N$ is smooth, there is a short exact sequence of cones on $Y$
$$\xymatrix@1{0\rto & T_{N/M}|_Y\rto & C_{Y/N}\rto & C_{X/M}|_Y\rto &
  0}\,.$$
This shows that $\c_{X/M}|_Y=\c_{Y/N}$. 

Comparing any two embeddings of $X$ with the diagonal, we get from
this the uniqueness of $\c_X$ for embeddable $X$.  Then we deduce that
$\c_X$ commutes with \'etale maps and thus glues with respect to the
\'etale topology. 
\end{pf}

If $X$ is smooth, $\c_X=(-1)^{\dim X}[X]$.

\begin{prop}
The basic properties of $\c_X$ are as follows:

(i) if $f:X\to Y$ is a smooth morphism of Deligne-Mumford stacks, then
$f^\ast\c_Y=(-1)^{\dim X/Y}\c_X$. 

(ii) if $X$ and $Y$ are Deligne-Mumford stacks, then $\c_{X\times
  Y}=\c_X\times\c_Y$. 
\end{prop}
\begin{pf}
Both of these follow from the product property of normal cones:
$C_{X/M}\times C_{Y/N}= C_{X\times Y/M\times N}$.
\end{pf}

\begin{rmk}
Maybe it would be appropriate to call $\c_X$ the {\em distinguished
  cycle }of $X$, in view of its relation to distinguished varieties
  in intersection theory (Definition~6.1.2 in~\cite{fulton}).
\end{rmk}

\subsection{The Euler obstruction $\nu_X$ of $\c_X$}

Consider MacPherson's local Euler obstruction $\Eu:Z_\ast(X)\to\Con(X)$,
which maps integral algebraic cycles on $X$ to constructible integer-valued
functions on 
$X$. Because $Eu$ commutes with \'etale maps and both $Z_\ast$ and
$\Con$ are sheaves with respect to the \'etale topology, $Eu$ is
well-defined for Deligne-Mumford stacks $X$ and defines an isomorphism
$Z_\ast(X)\to\Con(X)$. 

If $V$ is a prime cycle of dimension $p$ on the Deligne-Mumford stack
$X$, the constructible function $\Eu(V)$ takes the value
\begin{equation}\Label{leo}
\int_{\mu^{-1}(P)}c(\widetilde{T})\cap s\big(\mu^{-1}(P),\widetilde
V\big)
\end{equation}
 at the point $P\in X$. Here $\mu:\widetilde V\to V$ is the
Nash blowup of $V$ (the unique integral closed substack dominating $V$
of the Grassmannian of rank $p$ quotients of $\Omega_V$) and
$\widetilde T$ is the dual of the universal quotient bundle.  Moreover,
$c$ is the total Chern class and $s$ the Segre class of the normal
cone to a closed immersion. A proof that $\Eu(V)$ is constructible can
be found in~\cite{kenn}.

\begin{defn}
Let $X$ be a Deligne-Mumford stack. We introduce the canonical
integer valued constructible function
$$\nu_X=\Eu(\c_X)$$ 
on $X$.
\end{defn}

On a smooth stack $X$, the function $\nu_X$ is
locally constant, equal to $(-1)^{\dim X}$.

\begin{prop}
Let $X$ and $Y$ be Deligne-Mumford stacks.

(i) if $f:X\to Y$ is a smooth morphism, 
then $f^\ast \nu_Y=(-1)^{\dim X/Y} \nu_X$,

(ii) $\nu_{X\times Y}=\nu_X\boxdot \nu_Y$.
\end{prop}
\begin{pf}
Both facts follow from the compatibility of $\Eu$ with products.
\end{pf}

\subsubsection{Relation with Milnor numbers and vanishing
  cycles}

Suppose $M$ is a smooth scheme and $f:M\to\A^1$ is a regular function. Let
$X=Z(df)$ be the critical locus of $f$.  Then for any $\cc$-valued
point $P$ of $X$, we have
\begin{equation}\Label{fp}
\nu_X(P)=(-1)^{\dim M}\big(1-\chi(F_P)\big)\,,
\end{equation}
where $\chi(F_P)$ is the Euler characteristic of the Milnor fibre of
$f$ at
$P$.  For the proof, see~\cite{pp}, Corollary~2.4~(iii). Hence our
constructible function $\nu_X$ is equal to the function denoted $\mu$
in the literature (see, for example, [ibid.]).

Let $\Phi_f$ be the perverse sheaf of vanishing cycles on $X$.
It is a constructible complex on $X$ and therefore has a fibrewise
Euler characteristic 
$$\chi(\Phi_f)(P)=\sum_i (-1)^i\dim H^i_{\{P\}}(X,\Phi_f)\,,$$
which is a constructible function
on $X$. As a consequence of~(\ref{fp}), we have
\begin{equation}\Label{phi}
\nu_X=(-1)^{\dim M-1}\chi(\Phi_f)\,.
\end{equation}

\subsection{Weighted Euler characteristics}\Label{wec}

The Euler characteristic with compact supports $\chi$ is a 
$\zz$-valued 
function on   isomorphism classes of pairs $(X,f)$, where
$X$ is a scheme and $f$ a constructible function on $X$. It satisfies
the properties

(i) if $X$ is separated and smooth, $\chi(X,1)=\chi(X)$ is the usual
topological Euler characteristic of $X$,

(ii) $\chi(X,f+g)=\chi(X,f)+\chi(X,g)$,

(iii) if $X$ is the disjoint union of a closed subscheme $Z$ and its
open complement $U$, then $\chi(X,f)=\chi(U,f|_U)+\chi(Z,f|_Z)$,

(iv) $\chi(X\times Y,f\boxdot g)=\chi(X,f)\,\chi(Y,g)$,

(v) if $X\to Y$ is a finite \'etale morphism of degree $d$, then
$\chi(X,f|_X)=d\,\chi(Y,f)$, for any constructible function $f$ on
$Y$.  

These properties suffice to prove that $\chi$ extends uniquely to a
$\qq$-valued function on pairs $(X,f)$, where $X$ is a Deligne-Mumford
stack and $f$ a $\zz$-valued constructible function on $X$. (Use the
fact that every Deligne-Mumford stack is generically the quotient of a
scheme by a finite group, Corollaire~6.1.1 in~\cite{lmb}.) Properties
(i),\ldots,(v) continue to hold. 
We write $\chi(X)$ for $\chi(X,1)$. The rational number $\chi(X)$ is
often called the {\em orbifold Euler characteristic }of $X$.

\begin{prop}[Gau\ss-Bonnet]
If the
Deligne-Mumford stack $X$ is smooth and proper, we have that
$$\chi(X)=\int_{[X]} e(T_X)\,,$$
the integral of the Euler class (top Chern class) of the tangent
bundle. 
\end{prop}
\begin{pf}
First one proves that $\chi(I_X)=\chi(\ol{X})$, where $I_X$ is the
inertia stack and $\ol X$ the coarse moduli space of $X$. This can be
done by passing to stratifications and is therefore not difficult.
Then, invoking the Lefschetz trace formula for the identify on a
smooth and proper $X$ we get (see~\cite{trieste} for details)
$$\int_{[I_X]}e(T_{I_X})=\sum_i(-1)^i\dim H^i(X,\cc)=\chi(\ol X)\,.$$
Putting these two remarks together, we get the Gau\ss-Bonnet theorem
for $I_X$. To prove the theorem for $X$, note that the part of $I_X$
whose dimension is equal to $\dim X$ is a closed
and open substack $Y$ of $I_X$, which comes with a finite \'etale 
representable morphism $Y\to X$.  By induction on the dimension, the
theorem holds for $Y$. Then it also holds for $X$ by Property~(v) of
the Euler characteristic.
\end{pf}

Note that for stacks, $\chi(X)$ differs from the alternating
sum of the dimensions of the compactly supported cohomology groups.

\begin{defn}
Let $X$ be a Deligne-Mumford stack. 
Introduce the weighted Euler characteristic 
$$\tilde\chi(X)=\chi(X,\nu_X)=\chi\big(X,\Eu(\c_X)\big)\in\qq\,.$$
More generally, given a morphism $Z\to X$, define
$$\tilde\chi(Z,X)=\chi(Z,\nu_X|_Z)\,.$$
This definition is particularly useful for locally closed substacks
$Z\subset X$. 
\end{defn}

\begin{prop}
The weighted Euler characteristic $\tilde\chi(Z,X)$ has the basic properties:

(i) if $X$ is smooth, $\tilde\chi(Z,X)=(-1)^{\dim X}\chi(Z)$,

(ii) if $Z\to X$ is smooth, $\tilde\chi(Z,X)=(-1)^{\dim Z/X}\tilde\chi(Z)$. 

(iii) if $Z=Z_1\cup Z_2$ is the disjoint union of two locally closed
substacks, $\tilde\chi(Z_1,X)+\tilde\chi(Z_2,X)=\tilde\chi(Z,X)$, 

(iv) $\tilde\chi(Z_1\times Z_2,X_1\times
X_2)=\tilde\chi(Z_1,X_1)\,\tilde\chi(Z_2,X_2)$,

(v) given a commutative diagram
$$\xymatrix{Z\dto\rto & X\dto\\ W\rto & Y}$$
with $X\to Y$ smooth and $Z\to W$ finite \'etale, we have that
$\tilde\chi(Z,X)=(-1)^{\dim X/Y}\deg(Z/W)\,\tilde\chi(W,Y)$. 
\end{prop}
\begin{pf}
All these properties follow by combining the properties of $\nu_X$
with those of the orbifold Euler characteristic. 
\end{pf}

\begin{rmk}
Suppose $X$ is the disjoint union of an open substack $U$ and a closed
substack $Z$.  We have
$\tilde\chi(X)=\tilde\chi(U)+\tilde\chi(Z,X)$.  But because, in
general, $\tilde\chi(Z,X)\not=\tilde\chi(Z)$, we also have
$\tilde\chi(X)\not=\tilde\chi(U)+\tilde\chi(Z)$. 
\end{rmk}

\begin{rmk}
If $X$ is smooth and proper, 
$$\tilde\chi(X)=\int_{[X]}e(\Omega_X)\,.$$
Both Proposition~\ref{macaluprop} and Theorem~\ref{finally} can be
viewed as generalizations of this formula.
\end{rmk}

\subsection{The  Aluffi class}

The Mather class is a homomorphism $c^M:Z_\ast(X)\to A_\ast(X)$. It exists
for Deligne-Mumford stacks as well as for schemes.  The definition is
a globalization of the construction of the local Euler obstruction.
For a prime cycle $V$ of degree $p$ on $X$, we have
$$c^M(V)=\mu_\ast\big(c(\widetilde T)\cap[\widetilde V]\big)\,,$$
with the same notation as in~(\ref{leo}). 
We will only  need to use  the degree
zero part $c^M_0:Z_\ast(X)\to A_0(X)$.

\begin{defn}
Applying $c^M$ to our cycle $\c_X$, we obtain the {\bf Aluffi class}
$$\alpha_X=c^M(\c_X)\in A_\ast(X)\,,$$
\end{defn}

The class $\alpha_X$ was introduced by Aluffi~\cite{Alu}, although one
should note that the sign conventions in [ibid.] differ from ours.

If $X$ is smooth, its Aluffi class equals 
$$\alpha_X=(-1)^{\dim
  X}c(T_X)\cap[X]=c(\Omega_X)\cap[X]\,.$$ 

\begin{prop}\Label{macaluprop}
Let $X$ be a proper Deligne-Mumford stack. The formula (note that only
the degree zero 
component of $\alpha_X$ enters into it)
\begin{equation}\Label{macalu}
\tilde\chi(X)=\int_X\alpha_X\,
\end{equation}
is true in the following cases:

(i) if $X$ is a global finite group quotient,

(ii) if $X$ is a gerbe over a scheme,

(iii) if $X$ is smooth.
\end{prop}
\begin{pf}
In the smooth case, Formula~(\ref{macalu}) is the Gau\ss-Bonnet
theorem. 

For schemes, the proposition is true by a direct application of
MacPherson's theorem, which says that 
$$\chi\big(X,\Eu(\c)\big)=\int_X c^M(\c)\,,$$
for any cycle $\c\in Z_\ast(X)$. 

Let $f:X\to Y$ be a finite \'etale morphism of Deligne-Mumford
stacks, representable or not.  Then both the local Euler obstruction
and the Chern-Mather class commute with pulling back via $f$. It
follows that $\tilde\chi(X)=d\,\tilde\chi(Y)$ and
$\int_X\alpha_X=d\int_Y\alpha_Y$, where $d\in\qq$ is the degree of
$f$. Thus Formula~(\ref{macalu}) holds for $X$ if and only if it holds
for $Y$. 

If $X=[Y/G]$ is a global quotient of a scheme $Y$ by a finite group
$G$, there is the finite \'etale morphism $Y\to X$, proving Case~(i).
If $X$ is a gerbe, the morphism $X\to\ol X$ from $X$ to its coarse
moduli space is finite \'etale, proving Case~(ii).
\end{pf}

\begin{rmk}
Of course, it is very tempting to conjecture Formula~(\ref{macalu})
to hold true in general. 
\end{rmk}

\eject
\section{Remarks  on virtual cycle classes}

Let $X$  denote a scheme or a Deligne-Mumford stack. Let $L_X$ be the
cotangent complex of $X$. Recall from \cite{BF} that a {\em perfect
  obstruction theory }for $X$ is a derived category 
morphism $\phi:E\to L_X$, such that

(i) $E\in D(\O_X)$ is perfect, of perfect amplitude contained in the
interval $[-1,0]$, 

(ii) $\phi$ induces an isomorphism on $h^0$ and an epimorphism on
$h^{-1}$. 

\noindent Let us fix a perfect obstruction theory $E\to L_X$ for $X$.

Recall that $E$ defines a vector bundle stack $\EE$ over $X$: whenever
we write $E$ locally as a complex of vector bundles $E=[E_1\to E_0]$,
the stack $\EE$ becomes the stack quotient $\EE=[E_1^\vee/E_0^\vee]$. 

Recall also the {\em intrinsic normal cone }$\CC_X$.  Whenever $U\to
X$ is \'etale and $U\to M$ a closed immersion into a smooth scheme $M$,
the pullback $\CC_X|_U$ is canonically isomorphic to the stack
quotient $[C_{U/M}/(T_M|_U)]$, where $C_{U/M}$ is the normal cone. The
morphism $E\to L_X$ defines a closed immersion of cone stacks
$\CC_X\hookrightarrow\EE$. 

Recall, finally, that the obstruction theory $E\to L_X$ defines a {\em
  virtual fundamental class }$[X]^\mathrm{vir}\in A_{\rk E}(X)$, as the
  intersection of the fundamental class $[\CC_X]$ with the zero
  section of $\EE$: $$[X]^\mathrm{vir}=0^!_{\EE}[\CC_X]\,.$$
(For the last statement in the absence of global resolutions,
  see~\cite{Kresch}.)  Here $A_r(X)$ denotes the Chow group of
  $r$-cycles modulo rational equivalence on $X$ with
  values in $\zz$.

\subsection{Obstruction cones}

\begin{defn} We call $ob=h^1(E^\vee)$ the {\bf  obstruction sheaf }of
  the obstruction theory $E\to L_X$. 
\end{defn}

A {\em local resolution }of $E$ is a derived category homomorphism
$F\to E^\vee[1]|_U$, over some \'etale open subset $U$ of $X$, where $F$ is
a vector bundle over $U$ and the homomorphism $F\to E^\vee[1]|_U$ is such
that its cone is a locally free sheaf over $U$ concentrated in degree
$-1$. Alternatively, a local resolution may be defined as a local
presentation $F\to\EE|_U$ (over an \'etale open $U$ of $X$) of the vector
bundle stack $\EE$ associated to $E$. 

Recall that for every local resolution $F\to E^\vee[1]|_U$ there is an
associated cone $C\hookrightarrow F$, the {\em obstruction cone},
defined via the cartesian diagram of cone stacks over $U$
$$\xymatrix{
C\dto\rto\ar@{}[dr]|{\Box} & F\dto\\
\CC|_U\rto & \EE|_U}$$
where $\CC$ is the intrinsic normal cone of $X$.

Note that every local resolution $F\to E^\vee[1]|_U$ comes with a canonical
epimorphism of coherent sheaves $F\to ob|_U$.

\begin{prop}
let $\Omega$ be a vector bundle on $X$ and $\Omega\to ob$ an
epimorphism of coherent sheaves. Then there exists a unique closed
subcone $C\subset\Omega$ such 
that for every local resolution $F\to E^\vee[1]|_U$, with obstruction cone
$C'\subset F$,  and every lift $\phi$ 
\begin{equation}\Label{lift}
\vcenter{\xymatrix{
& F\dto\\
\Omega|_U\rto\urto^\phi & ob|_U}}
\end{equation}
we have that $C|_U=\phi^{-1}(C')$, in the scheme-theoretic sense.
\end{prop}
\begin{pf}
\'Etale locally on $X$, presentations $F$ and lifts $\phi$ always
exist. The uniqueness of $C$ follows. 

So far, we have only considered $ob$ as a coherent sheaf on $X$.  We
can extend it to a sheaf on the big \'etale site of $X$ in the
canonical way. We may then think of $ob$ as the coarse moduli sheaf of
$\EE$. Let $cv$ be the coarse moduli sheaf of the intrinsic normal
cone $\CC$.  The key facts are

(i) $cv\hookrightarrow ob$ is a subsheaf,

(ii) the diagram
\begin{equation}\Label{cst}
\vcenter{\xymatrix{
\CC\dto\ar@{^{(}->}[r]\ar@{}[dr]|{\Box} & \EE\dto\\
cv\ar@{^{(}->}[r] & ob}}
\end{equation}
is a cartesian diagram of stacks over $X$.

Both of these facts are local in the \'etale topology of $X$, so we
may assume that $E$ has a global resolution $E^\vee=[H\to F]$.  Let
$C'\subset F$ be the obstruction cone.  Then $C'$ is invariant under the
action of $H$ on $F$. Note that $ob$ is the sheaf-theoretic quotient
of $F$ by $H$ and $cv$ the sheaf-theoretic quotient of $C'$ by $H$.
Simple sheaf theory on the big \'etale site of $X$ (exactness of the
associated sheaf functor) implies 

(i) $cv\hookrightarrow ob$ is a subsheaf,

(ii)  the diagram
\begin{equation}\Label{first}
\vcenter{\xymatrix{
C'\dto\ar@{^{(}->}[r]\ar@{}[dr]|{\Box} & F\dto\\
cv\ar@{^{(}->}[r] & ob}}
\end{equation}
is a cartesian diagram of sheaves on the big \'etale site of $X$.
This implies that 
Diagram~(\ref{cst}) is cartesian, proving the key facts.

We now construct the subsheaf $C\subset \Omega$ as the fibred product of
sheaves on the big \'etale site of $X$
\begin{equation}\Label{second}
\vcenter{\xymatrix{
C\rto\dto\ar@{}[dr]|{\Box} & \Omega\dto\\
cv\rto & ob}}
\end{equation}
Then any diagram such as~(\ref{lift}) gives rise to a cartesian
diagram of big \'etale sheaves
$$\xymatrix{
C\rto \dto\ar@{}[dr]|{\Box} & \Omega\dto^\phi\\
C'\rto & F}$$
This latter diagram is cartesian, because Diagrams~(\ref{first})
and~(\ref{second}) are.  This proves the claimed property of $C$, as
well as the fact that $C$ is a closed subcone of $\Omega$, in the
scheme-theoretic sense. 
\end{pf}

\begin{rmk}
In \cite{BF}, it was shown that the subsheaf $cv\hookrightarrow ob$
classifies small curvilinear obstructions.
\end{rmk}

\begin{defn}\Label{defnoc}
We call $C\subset\Omega$ the {\bf obstruction cone }associated to the
epimorphism $\Omega\to ob$. 
\end{defn}

Note that $ob$ is in general bigger than the actual sheaf of obstructions,
which is the abelian subsheaf of $ob$ generated by $cv$.  Thus $cv$ is
intrinsic to $X$, whereas $ob$ depends on $E\to L_X$.

\subsection{The virtual fundamental class}

\newcommand{\Qcoh}{\text{\rm Qcoh-$\O_X$}}
\newcommand{\qcoh}{\rm qcoh}
\newcommand{\coh}{\rm coh}

\begin{lem}
If $X$ is a quasi-projective Deligne-Mumford stack, every perfect
obstruction theory $E\to L_X$ has a global resolution.
\end{lem}
\begin{pf}
Let $D(\O_X)$ be the derived category of  sheaves of
$\O_X$-modules on the (small) \'etale site of $X$ and let $D(\Qcoh)$
be the derived category of the category of quasi-coherent $\O_X$-modules. 
First prove that the natural functor $D^+(\Qcoh)\to D^+_{\qcoh}(\O_X)$
is an equivalence of categories.  For this, show that the inclusion of
categories $(\Qcoh)\to(\text{\rm $\O_X$-mods})$ has a right adjoint
$Q$, which commutes with pushforward along morphisms of
quasi-projective stacks. Prove that quasi-coherent sheaves are acyclic
for $Q$ and satisfy that $QF\to F$ is an isomorphism. Thus the
right derivation of $Q$ provides a quasi-inverse to the 
inclusion.  To reduce all these claims to the affine case use a
groupoid $U_1\toto U_0$  presenting $X$, which is \'etale and has
affine $U_1$, $U_0$.  For the details of the proof, see Section~3 of
Expos\'e~II in SGA6.

Next, prove that every quasi-coherent sheaf on $X$ is a direct limit of
coherent sheaves.  For this, it is convenient to choose a finite flat
cover $f:Y\to X$, where $Y$ is a quasi-projective scheme. Construct a
right adjoint $f^!$ to $f_\ast$, from $(\Qcoh)$ to $({\text{\rm
    Qcoh-$\O_Y$}})$. (This can be done \'etale locally over $X$.) 
Now, let $F$ be a quasi-coherent $\O_X$-module. There exist coherent
sheaves $G_i$ on $Y$, such that $f^! F=\lim_{\to} G_i$, because $Y$ is
quasi-projective. Since it
admits a right adjoint, $f_\ast$ commutes with direct limits and so we
have $f_\ast f^! F=\lim_{\to} f_\ast G_i$. The trace map $f_\ast f^!
F\to F$ is onto, and so we get a surjection $\lim_{\to}f_\ast G_i\to
F$.  Since the $f_\ast G_i$ are coherent, $F$ is, indeed, an
inductive limit of coherent modules. 

With these preparations, we can now construct a global resolution of
the perfect complex $E\in D^b_{\rm coh}(\O_X)$.  First, we may assume
that $E$ is given by a 2-term complex $E=[E_{1}\to E_0]$, where $E_0$
and $E_1$ are quasi-coherent. (Our above argument gives an infinite
complex of quasi-coherents, which we may cut off, because the kernel
of a morphism between quasi-coherent sheaves is quasi-coherent.)  Then
we choose coherent sheaves $G_i$ on $X$ such that
$E_0=\lim_{\to}G_i$. The images of the $G_i$ in $h^0(E)=\cok(E_{1}\to
E_0)$ stabilize, because $X$ is noetherian and hence the coherent
$\O_X$-module $h^0(E)$ satisfies the ascending chain condition. So
there exists a coherent sheaf $G_i\to E_0$, which maps surjectively
onto $h^0(E)$. Now find a locally free coherent $F_0$ mapping onto
$G_i$ (and hence onto $h^0(E)$), and define $F_{1}=E_{1}\times_{E_0}
F_0$.  Then $F_{1}$ is 
automatically locally free coherent and $F=[F_{1}\to F_0]$ maps
quasi-isomorphically to $[E_{1}\to E_0]$.  Thus $F$ provides us with
the required global resolution of $E$.
\end{pf}

\begin{prop}\Label{vireq}
Consider a quasi-projective Deligne-Mumford stack $X$ and a perfect
obstruction theory $E\to L_X$ with obstruction sheaf $ob$.  Let
$\Omega$ be a vector bundle over $X$ and $\Omega\to ob$ an epimorphism
of coherent sheaves.  Let $C\subset\Omega$ be the associated
obstruction cone.  Then $C$ is of pure dimension $\rk E +\rk \Omega$
and we have
$$[X]^\mathrm{vir}=0_\Omega^![C]\,.$$
\end{prop}
\begin{pf}
Let $F\to E^\vee[1]$ be a global resolution of $E$  with obstruction
cone $C'\subset F$. 
Start by constructing the fibred product of coherent sheaves
\begin{equation}\Label{small}
\vcenter{\xymatrix{
\mathcal{P}\dto\rto\ar@{}[dr]|{\Box} & F\dto\\
\Omega\rto & ob}}
\end{equation}
and choosing an epimorphism of coherent sheaves
$F'\to\mathcal{P}$, where $F'$ is locally free. 
Wet get a commutative diagram of sheaf epimorphisms
$$\xymatrix{
F'\ar@{->>}[d]\ar@{->>}[r] & F\ar@{->>}[d]\\
\Omega\ar@{->>}[r] & ob}$$
which we can now consider as a diagram of sheaves on the big \'etale
site of $X$.

{\em Remark.}  Diagram (\ref{small}) is a cartesian diagram of sheaves
on the small \'etale site of $X$.  This is because fibred products of
(small) coherent sheaves do not commute with base change, and so if we
had taken the fibred product of big sheaves, $\mathcal{P}$ would not
have ended up coherent.  After having chosen $F'$, we do not any
longer have use for the cartesian property of the diagram, and so we
pass back to big sheaves, as commutativity of diagrams and the
property of being an epimorphism are stable under base change.

Now, of course, $F'\to \Omega$ and $F'\to F$ are epimorphisms of
vector bundles. The preimage of $cv\hookrightarrow ob$ in $\Omega$ is
$C$, and in $F$ is $C'$.  It follows that $C'$ and $C$ have the
same preimage in $F'$.  This implies by standard arguments the claim
about the dimension of $C$ and the fact that
$[C']\cap[O_F]= [C]\cap[O_\Omega]$. 
\end{pf}

\eject
\section{Symmetric obstruction theories}

We will summarize the main features of symmetric obstruction
theories.  For proofs, see~\cite{bf5}. 
Throughout this section, $X$ will denote a Deligne-Mumford stack.

\subsection{Non-degenerate symmetric bilinear forms}

\begin{defn}
Let $E\in
D^b_{\coh}(\O_X)$ be a perfect 
complex.  A {\bf non-degenerate symmetric bilinear form of degree 1 }on $E$ is
an isomorphism $\theta:E\to E^\vee[1]$, satisfying
$\theta^\vee[1]=\theta$.
\end{defn}

Of course, it has to be understood that $\theta$ is a morphism in the
derived category, and invertible as such.  The duals appearing in the
definition are derived.

\begin{ex}\Label{simple}
A simple example of a perfect complex with non-degenerate symmetric
bilinear form of degree 1 is given as follows. Let $F$ be a vector
bundle on $X$, endowed with a symmetric bilinear form, inducing a
homomorphism $\alpha:F\to F^\vee$. To define the complex $E=[ F\to
F^\vee]$, put $F^\vee$ in degree $0$ and $F$ in degree $-1$. Since the
components of $E$ are locally free, we can compute the derived dual as
componentwise dual. We fine $E^\vee[1]=E$. So we may and will define
$\theta$ to be the identity, i.e., $\theta_{1}=\id_{F}$ and
$\theta_0=\id_{F^\vee}$:
$$\xymatrix{
E\dto^\theta \ar@{}[r]|{=} & [F\rto^\alpha\dto_{1} & F^\vee\dto^1]\\
E^\vee[1] \ar@{}[r]|{=} & [F\rto^{\alpha} & F^\vee]}$$
Note that $\theta$ is an
isomorphism, whether or not $\alpha$ is non-degenerate.
\end{ex}

\begin{ex}
As a special case of Example~\ref{simple}, consider a regular function
$f$ on a smooth variety 
$M$. The Hessian of $f$ defines a symmetric bilinear form on
$T_M|_X$, where $X=Z(df)$. Hence we get a non-degenerate symmetric
bilinear form
on the complex $[T_M|_X\to\Omega_M|_X]$, which is, by the way, a perfect
obstruction theory for $X$.
\end{ex}

\begin{defn}
Let $A$ and $B$ be perfect complexes endowed with non-degenerate
symmetric forms $\theta:A\to A^\vee[1]$ and $\eta:B\to B^\vee[1]$. An
{\bf   isometry }$\Phi:(B,\eta)\to(A,\theta)$ is an isomorphism
$\Phi:B\to A$, such that the diagram
$$\xymatrix{
 B\dto_{\eta} \rrto^\Phi && A\dto^\theta\\
 B^\vee[1] && A^\vee[1]\llto_{\Phi^\vee[1]}}$$
commutes in $D(\O_X)$. Since $\eta$ and $\theta$ are isomorphisms, this
amounts to saying that $\Phi^{-1}=\Phi^\vee[1]$ (if we use
$\eta$ and $\theta$ to identify).
\end{defn}

\subsection{Symmetric obstruction theories}

\begin{defn}
A perfect obstruction theory $E\to L_X$ for $X$ is called {\bf symmetric},
if $E$ is endowed with a non-degenerate symmetric bilinear form
$\theta:E\to E^\vee[1]$. 
\end{defn}

If $E$ is symmetric, we have $\rk E=\rk (E^\vee[1])=-\rk
E^\vee=-\rk E$ and hence $\rk E=0$. So the expected dimension is zero.
Therefore, we can make the following definition:

\begin{defn}
Let $X$ be endowed with a symmetric obstruction theory and assume that
$X$ is proper. The {\bf virtual count }(or {\em Donaldson-Thomas type
  invariant}) of $X$ is the number
$$\#^\mathrm{vir}(X)=\deg [X]^\mathrm{vir}=\int_{[X]^\mathrm{vir}}1\,.$$
If $X$ is a scheme (or an algebraic space), $\#(X)$ is an integer,
otherwise a rational number.
\end{defn}

\begin{rmk}
For a symmetric obstruction theory $E\to L_X$, we have
$ob=h^1(E^\vee)=h^0(E^\vee[1])=h^0(E)=\Omega_X$. So 
the obstruction sheaf
is canonically isomorphic to the sheaf of differentials.
\end{rmk}

\begin{rmk}\Label{roc}
Let $X$ be endowed with a symmetric obstruction theory.  Then for any
closed embedding $X\to M$ into a smooth Deligne-Mumford stack $M$, we
get a canonical epimorphism of coherent sheaves
$\Omega_M|_X\to\Omega_X=ob$, and hence a canonical closed subcone
$C\hookrightarrow \Omega_M|_X$, the obstruction cone of
Definition~\ref{defnoc}. Via the inclusion
$\Omega_M|_X\hookrightarrow\Omega_M$ we think of $C$ as a closed conic
substack of $\Omega_M$. If $X$ is quasi-projective,
Proposition~\ref{vireq} applies and so $C$ is
pure dimensional and $\dim C=\dim M=\frac{1}{2}\dim\Omega_M$. We will
show below that $C$ is {\em Lagrangian}. 
\end{rmk}

\begin{rmk}
Any symmetric obstruction theory on $X$ induces (by restriction) in a
canonical way a symmetric obstruction theory on $U\to X$, for every
\'etale morphism $U\to X$. 
\end{rmk}

\begin{rmk}
If $E$ is a symmetric obstruction theory for $X$ and $F$ a symmetric
obstruction theory for $Y$, then $E\boxplus F$ (see \cite{BF}) is
naturally a symmetric obstruction theory for $X\times Y$.
\end{rmk}

\subsection{Examples}\Label{exam}

\subsubsection{Lagrangian intersections}

Let $M$ be a complex symplectic manifold and $V$, $W$ two Lagrangian
submanifolds.  Let $X$ be their scheme-theoretic intersection. Then $X$
carries a canonical symmetric obstruction theory. This generalizes to
Deligne-Mumford stacks.

\subsubsection{Sheaves on Calabi-Yau threefolds}

Let $Y$ be a smooth projective Calabi-Yau threefold and $L$ a line bundle
on $Y$.  Let $X$ be any open substack of the stack of stable sheaves
of positive rank $r$ with determinant $L$. (For example, $X$ could be the
stack of sheaves 
of a fixed Hilbert polynomial admitting no strictly semi-stable
sheaves.  Then $X$ would be proper.) 

Let $\mathcal{E}$ be the universal sheaf and $\mathcal{F}$ the shifted
cone of the trace map:
$$\xymatrix@C=1pc{
& \O\dlto_{+1} & \\
\mathcal{F}\rrto && R\sheafhom(\mathcal{E},\mathcal{E})\ulto_{tr}}$$
Then $R\pi_\ast \mathcal{F}[2]$ is a symmetric obstruction theory for
$X$.  Here $\pi:Y\times X\to X$ is the projection.  For the proof,
see~\cite{RT} or~\cite{bf5}.

Note that $X$ is a $\mu_r$-gerbe over a quasi-projective
scheme. Moreover, $X$ is a quasi-projective stack.

\subsubsection{Hilbert schemes of local Calabi-Yau threefolds}

If we restrict to rank one sheaves, we can consider the following more
general situation.  Let $Y$ be a smooth projective threefold with a
section of the 
anticanonical bundle whose zero locus we denote by $D$.  Recall that
stable sheaves with trivial determinant can be considered as ideal
sheaves on $Y$.  

Let $X$ be an open subscheme
of the Hilbert scheme of ideal sheaves on $Y$.  We require that $X$ 
consists entirely of ideal sheaves whose associated subschemes of $Y$
are disjoint from $D$. Then, with the same notation as above,
$R\pi_\ast\mathcal{F}[2]$ is a symmetric obstruction theory for $X$.

Note that $X$ is a quasi-projective scheme.

\subsubsection{Stable maps}

Let $Y$ be a Calabi-Yau threefold. Let $X$ be the open substack of the
stack of stable maps $\ol{M}_{g,n}(Y,\beta)$, corresponding to stable
maps which are  immersions from a smooth curve to $Y$.  Then 
the Gromov-Witten obstruction theory for $\ol{M}_{g,n}(Y,\beta)$ is
symmetric  over $X$.

\subsection{Local structure: almost closed 1-forms}

\begin{defn}
A differential form $\omega$ on a smooth Deligne-Mumford stack $M$ is
called {\bf almost closed}, if $d\omega\in I\Omega_M^2$. Here $I$ is the
ideal sheaf of the zero locus of $\omega$ (in other words the image of
$\omega^\vee:T_M\to \O_M$).  Equivalently, we may say that $d\omega|_X=0$
as a section of $\Omega_M^2|_X$, where
$X$ is the zero locus of $\omega$, i.e., $\O_X=\O_M/I$. 
\end{defn}

Of course, in local coordinates $x_1,\ldots,x_n$, where $\omega=\sum_i
f_i\,dx_i$, being almost closed means that 
$$\frac{\del f_i}{\del x_j}\equiv\frac{\del f_j}{\del x_i}\mod
(f_1,\ldots,f_n)\,,$$
for all $i,j=1,\ldots,n$.

\begin{rmk}
Let $M$ be a smooth Deligne-Mumford stack and $\omega$ an almost
closed 1-form on $M$ with zero locus $X=Z(\omega)$.
It is a general principle, that a section of a vector bundle defines a
perfect obstruction theory for the zero locus of the section. 
In our case, this obstruction theory is given by
$$\xymatrix{
E\dto\ar@{}[r]|{=}& [T_M|_X \rto^{d\circ\omega^\vee}\dto_{\omega^\vee}
& \Omega_M|_X]\dto^1\\ 
\tau_{\geq-1}L_X\ar@{}[r]|{=} & [I/I^2\rto^{d} & \Omega_M|_X]}$$
We have only displayed the cutoff at -1 of $E\to L_X$, but that is the
only part of the obstruction theory that intervenes in our discussion.

This obstruction theory is symmetric, in a canonical way, because
under our assumption that $\omega$ is almost closed we have that
$d\circ\omega^\vee$ is self-dual, as a homomorphism of vector bundles
over $X$.  (See Example~\ref{simple}.) We denote this symmetric
obstruction theory by 
$H(\omega)\to L_X$, where
$$\xymatrix@1{
H(\omega)=[T_M|_X\rto^-{d\circ\omega^\vee}&\Omega_M|_X]}\,.$$
\end{rmk}

We will show that, at least locally, every symmetric obstruction
theory is given in this way by an almost closed 1-form. 

\begin{prop}\Label{preex}
Suppose $E\to L_X$ is a symmetric obstruction theory for the
Deligne-Mumford stack $X$. Then \'etale locally in $X$
(Zariski-locally if $X$ is a scheme) there exists a closed immersion
$X\hookrightarrow M$ 
of $X$ into a smooth scheme $M$ and an almost closed 1-form $\omega$
on $M$ and an isometry $E\to H(\omega)$ such that the diagram
$$\xymatrix@C=1pc{E\drto\rrto && H(\omega)\dlto\\
& L_X & }$$
commutes in the derived category of $X$.
\end{prop}
\begin{pf}
Let $P$ be a $\cc$-valued point of $X$. By passing to an \'etale
neighbourhood of $P$, we may assume given a closed immersion
$X\hookrightarrow M$ into a smooth scheme $M$ of dimension $\dim
M=\dim \Omega_X|_P$.  Moreover, we may assume that $E=[E_1\to E_0]$ is
given by a homomorphism of vector bundles such that $\rk E_0=\rk
E_1=\dim M$ 
and $E\to L_X$ is given by a homomorphism of complexes 
$$\xymatrix{
E\dto\ar@{}[r]|{=} & [E_1\dto_{\phi_1}\rto^\alpha & E_0\dto^{\phi_0}]\\
\tau_{\geq-1}L_X \ar@{}[r]|{=} & [I/I^2\rto & \Omega_M|_X]}$$
Since $\phi_0$ is an isomorphism at $P$, by passing to a smaller
neighbourhood of $P$, we may assume that $\phi_0$ is an isomorphism and
use it to identify $E_0$ with  $\Omega_M|_X$. 

For the  symmetric form
$\theta:E\to E^\vee[1]$ let us use notation 
$\theta=(\theta_1,\theta_0)$. Then the equality of derived category
morphisms $\theta^\vee[1]=\theta$ implies
that, locally, $\theta^\vee[1]=(\theta_0^\vee, \theta_1^\vee)$ and
$\theta=(\theta_1, \theta_0)$ are homotopic.  So let $h:E_0\to
E_0^\vee$ be a homotopy: 
\begin{align*}
h\alpha & =\theta_1-\theta_0^\vee\\
\alpha^\vee h &=\theta_0-\theta_1^\vee\,.
\end{align*}
Now define 
\begin{align*}
\lambda_0&={\textstyle\frac{1}{2}}(\theta_0+\theta_1^\vee)\\
\lambda_1&={\textstyle\frac{1}{2}}(\theta_1+\theta_0^\vee)\,.
\end{align*}
One checks that $(\lambda_1,\lambda_0)$ is a homomorphism of complexes, and as
such, homotopic to $(\theta_1,\theta_0)$. Thus $(\lambda_1,\lambda_0)$
represents the 
derived category morphism $\theta$, and has the property that
$\lambda_1=\lambda_0^\vee$:
$$\xymatrix{
E\dto_\theta\ar@{}[r]|{=} & [E_1\dto_{\lambda^\vee}\rto^\alpha &
  \Omega_M|_X\dto^{\lambda}]\\ 
E^\vee[1] \ar@{}[r]|{=} & [T_M|_X\rto^{\alpha^\vee} & E_1^\vee]}$$
Since $\theta$ is a quasi-isomorphism, $\lambda$ is necessarily an
isomorphism at $P$, hence, without loss of generality, an
isomorphism.  Use $\lambda$ to identify. Then we have written our
obstruction theory as
$$\xymatrix{
E\dto\ar@{}[r]|{=} & [T_M|_X\dto_{\phi_1}\rto^\alpha &
  \Omega_M|_X\dto^{1}]\\ 
\tau_{\geq-1}L_X \ar@{}[r]|{=} & [I/I^2\rto & \Omega_M|_X]}$$
with $\alpha=\alpha^\vee$.  Lift $\phi_1:T_M|_X\to I/I^2$ in an
arbitrary fashion to a homomorphism $\omega^\vee:T_M\to I$, defining
an almost closed 1-form $\omega$, such that $E=H(\omega)$.
\end{pf}

We need a slight amplification of this proposition:

\begin{cor}\Label{exists}
Let $E$ be a symmetric obstruction theory for $X$ and let
$X\hookrightarrow M'$ be an embedding into a smooth Deligne-Mumford
stack $M'$. Then \'etale
locally in $M'$, there exists an
almost closed 1-form $\omega$  
on $M'$,
such that $X=Z(\omega)$ and $E\to L_X$ is  isometric to $H(\omega)\to
L_X$.
\end{cor}
\begin{pf}
Let $P\in X$. The proof of Proposition~\ref{preex} actually gives $M$
is an \'etale slice though $P$ in $M'$.  Then write $M'$ locally as a
product of the slice with a complement to the slice.
\end{pf}

\eject
\section{Microlocal geometry}

\subsection{Conic Lagrangians inside $\Omega_M$}

Let $M$ be a smooth scheme. The cotangent bundle $\Omega_M$ carries
the {\em tautological }1-form $\alpha\in\Omega^1(\Omega_M)$.  It is
the image of the identity under
$\pi^\ast\Omega^1_M\to\Omega^1_{\Omega_M}$, the pullback map for
1-forms under the projection $\pi:\Omega_M\to M$.
Its differential $d\alpha$ defines the tautological symplectic
structure on $\Omega_M$. 

Let $\theta$ be the vector field on
$\Omega_M$ which generates the $\cc^\ast$-action on the fibres of
$\Omega_M$.  It is the image of the identity under
$\pi^\ast\Omega_M\to T_{\Omega_M}$, the map which identifies elements
of the vector bundle $\Omega_M$ with vertical tangent vectors for the
projection $\pi$. 

The basic relation between these tensors is
$$\alpha=d\alpha(\theta,\,\cdot\,)\,.$$

Any local \'etale coordinate system  $x_1,\ldots,x_n$ on $M$ induces
the {\em canonical }coordinate system $x_1,\ldots,x_n,p_1,\ldots,p_n$
on $\Omega_M$.  In such  canonical coordinates we have $\alpha=\sum_i
p_idx_i$ and $\theta=\sum_i p_i\frac{\del}{\del p_i}$. 

Consider an  irreducible closed subset $C\subset\Omega_M$.
We call $C$ {\em conic}, if
$\theta$ is tangent to $C$ at the generic point of $C$. We call $C$
{\em Lagrangian}, if $\dim C=\dim M$ and $d\alpha$ vanishes at the
generic point of $C$.

\begin{lem}
The irreducible closed subset $C\subset\Omega_M$ is conic and
Lagrangian if and only if $\dim C=\dim M$ and $\alpha$ vanishes at the
generic point of $C$. 
\end{lem}
\begin{pf}
Suppose $C$ is Lagrangian. The basic relation shows that $\alpha$
vanishes at smooth  points of $C$ if and only if $\theta\in T_C^\perp=T_C$
at such points.
\end{pf}

If $V\subset M$ is an irreducible closed subset, the  closure in
$\Omega_M$ of the cotangent bundle to any smooth dense open subset
of $V$ is conic Lagrangian. This already describes all conic Lagrangians:

\begin{lem}\Label{lagr}
Let $C\subset\Omega_M$ be a closed irreducible subset.  Let $V=\pi(C)$
be its image in $M$ and let $N\subset \Omega_M$ be the closure of the
cotangent bundle of any smooth dense open subset of $V$. 

If $C$ is conic and Lagrangian then it is equal to $N$.
\end{lem}
\begin{pf}
(See also \cite{kenn}, for a coordinate free proof.)
Choose local coordinates $x_1,\ldots,x_n$ for $M$ around a smooth
point of $V$, in such a way that $V$ is cut out by the equations
$x_1=\ldots,x_k=0$. Then $dx_{k+1}\ldots,dx_n$ are linearly
independent at the generic point of $V$. By generic smoothness of the
projection $C\to V$, these forms stay linearly independent at the
generic point of $C$.  Since $\alpha$ restricts to $\sum_{i=k+1}^n
p_idx_i$ at the generic point of $C$, and $\alpha$ vanishes there, we
see that $p_{k+1},\ldots,p_n$ vanish at the generic point of $C$.
Thus $x_1,\ldots,x_k,p_{k+1},\ldots p_k$ vanish at the generic point
of $C$. 

On the other hand, $N$ is cut out generically by
$x_1,\ldots,x_k,p_{k+1},\ldots p_k$.  Thus we have proved that the
generic point of $C$ is contained in $N$. Then $C=N$ for dimension
reasons. 
\end{pf}

\begin{defn}
A closed subset of $\Omega_M$ is called {\em conic Lagrangian},
if every one of its irreducible components is conic and Lagrangian.

An algebraic cycle on $\Omega_M$ is {\em conic Lagrangian} if its
support is conic Lagrangian, a closed subscheme 
of $\Omega_M$ is {\em conic Lagrangian }if its underlying closed
subset is conic Lagrangian. 
\end{defn}

\begin{rmk}
The property of being a conic Lagrangian is local in the \'etale topology of
$M$, so it makes sense also in the case when $M$ is a smooth
Deligne-Mumford stack. 
\end{rmk}

\subsubsection{Cycles}

Consider a smooth Deligne-Mumford stack $M$ of dimension $n$. Let
$\L(\Omega_M)\subset Z_n(\Omega_M)$ be the subgroup 
generated by the conic Lagrangian prime cycles.  
 
If $V$
is a prime cycle (integral closed substack) of $M$, we consider the
closure in $\Omega_M$ of the conormal bundle of any smooth dense open
subset of $V$ and denote it by $\ell(V)$.  Note that $\ell(V)$ is a conic
Lagrangian prime cycle on $\Omega_M$. 
This defines the homomorphism 
\begin{align}
\Label{LL}L:Z_\ast(M)&\longrightarrow \L(\Omega_M)\\
V&\longmapsto (-1)^{\dim V}\ell(V)\,.\nonumber
\end{align}

Conversely, if $W$ is a conic prime cycle on $\Omega_M$, intersecting
(set-theoretically) with the zero section of $\pi:\Omega_M\to M$ or
taking the (set-theoretic) image $\pi(W)$, we obtain the same prime
divisor in $M$.  Restricting to conic Lagrangian cycles, we obtain a
homomorphism
\begin{align}
\Label{PI}\pi:\L(\Omega_M)&\longrightarrow Z_\ast(M)\\
W&\longmapsto (-1)^{\dim\pi(W)}\pi(W)\,.\nonumber
\end{align}

By Lemma~\ref{lagr} the homomorphisms $L$ and $\pi$ between
$Z_\ast(M)$ and $\L(\Omega_M)$ are inverses of of each other.

\begin{rmk}
The characteristic cycle map $\Ch:\Con(X)\to\L(\Omega_M)$ is the
unique homomorphism making the diagram
$$\xymatrix{Z_\ast(M)\drto_\Eu\rrto^L && \L(\Omega_M)\\
&\Con(M)\urto_\Ch}$$
commute.  
\end{rmk}  

\subsubsection{Microlocal view of the Mather class}

\begin{prop}
Let $M$ be a smooth Deligne-Mumford stack.
The diagram
$$\xymatrix{
Z_\ast(M)\rto^L \drto_{c^M_0}& \L(\Omega_M)\dto^{0^!}\\
& A_0(M)}$$
commutes.
\end{prop}
\begin{pf}
Assume $V\subset M$ is a prime cycle of dimension $p$. Let $\mu:\widetilde
M\to M$ be the Grassmannian of rank-$p$ quotients of $\Omega_M$ and
$\nu:\widetilde V\to V$ the closure inside $\widetilde M$ of the canonical
rational section $V\dotarrow \widetilde M$. Then
$c^M_0(V)=(-1)^p\nu_\ast\big(c_p(Q)\cap [\widetilde
  V]\big)$, where $Q$ is the universal quotient bundle on $\widetilde
M$. 

Let us denote the kernel of the universal quotient map by $N$.  Then
on $\widetilde V$ we have the exact sequence of vector bundles
$$\xymatrix@1{0\rto & N|_{\widetilde V}\rto &
  \mu^\ast\Omega_M|_{\widetilde V}\rto & Q|_{\widetilde V}\rto  &
  0}\,.$$
It implies that $c_p(Q)\cap[\widetilde
  V]=0^!_{\Omega_M}[N|_{\widetilde V}]\in A_0(\widetilde V)$.
$$\xymatrix{
N|_{\widetilde V}\dto\rto& \mu^\ast\Omega_M\dto\rto&\Omega_M\dto &
C\lto\dto\\
\widetilde V\rto & \widetilde M\rto^\mu & M & V\lto}$$

Let $C=\ell(V)\subset\Omega_M$ be the conic Lagrangian prime cycle
defined by $V$. There is a canonical rational section
$C\dotarrow\mu^\ast\Omega_M$ and the closure of the image is equal to
$N|_{\widetilde V}$. Hence we have a projection map
$\eta:N|_{\widetilde V}\to C$ which is a proper birational map of
integral stacks. It fits into the cartesian diagram:
$$\xymatrix{
\widetilde V\dto_0\rto^\nu & V\rto\dto_0 & M\dto^0\\
N|_{\widetilde V}\rto^\eta & C\rto &\Omega_M}$$
Thus we have $\nu_\ast(0^!_{\Omega_M}[N|_{\widetilde
    V}])=0^!_{\Omega_M}\eta_\ast [N|_{\widetilde
    V}]=0^!_{\Omega_M}[C]$ and $c_0^M(V)=(-1)^p 0^!_{\Omega_M}[C]=
(-1)^p 0^!_{\Omega_M}\ell(V)= 0^!_{\Omega_M}L(V)$.
\end{pf}

\begin{cor}\Label{loc}
If $X$ is a closed substack of $M$, the diagram
$$\xymatrix{
Z_\ast(X)\rto^L \drto_{c^M_0}& \L_X(\Omega_M)\dto^{0_{\Omega_M}^!}\\
& A_0(X)}$$
commutes as well.  Here $\L_X(\Omega_M)$ denotes the subgroup of conic
Lagrangian cycles lying over cycles contained in $X$.
\end{cor}
\begin{pf}
We just have to remark that the Mather class computed inside $M$
agrees with the Mather class computed inside $X$.
\end{pf}

\begin{rmk}
This proves the existence of Diagram~(\ref{pretty}).
\end{rmk}

\subsection{The fundamental lemma on almost closed 1-forms}

Let $M$  denote a smooth scheme. 
Let $\omega$ be a 1-form on $M$ and $X=Z(\omega)$
its scheme-theoretic zero-locus. Considering $\omega$ as a linear
homomorphism $T_M\to \O_M$, its image $I\subset \O_M$ is the ideal
sheaf of $X$. The epimorphism $\omega^\vee:T_M\to I$ restricts to an
epimorphism $\omega^\vee:T_M|_X\to I/I^2$, which gives rise to a closed
immersion of cones $C_{X/M}\hookrightarrow \Omega_M|_X$. Via
$\Omega_M|_X\hookrightarrow \Omega_M$ we consider  $C=C_{X/M}$ as a subscheme
of $\Omega_M$.

\begin{thm}\Label{main}
If the 1-form $\omega$ is almost closed, the closed subscheme
$C\subset\Omega_M$ it defines is conic Lagrangian.
\end{thm}
The proof will follow after an example.

\begin{ex}
The case where the zero locus $X$ of $\omega$ is smooth is easy:
if $\omega$ is almost closed with smooth zero locus, $C\subset
\Omega_M$ is equal to the the conormal bundle $N_{X/M}^\vee\subset
\Omega_M$ and is hence conic Lagrangian. 

For the general case, this implies that all components of $C$ which
lie over smooth points of $X$ 
are conic Lagrangian. 
\end{ex}

\subsubsection{The proof of Theorem~\ref{main}}

We start with two lemmas. 

\begin{lem}\Label{value}
Let $B$ be an integral noetherian $\cc$-algebra, $f\in B$ non-zero and
$Q:B\to K$ a morphism to a field, such that $Q(f)=0$.  Then there
exists a field extension $L/K$, a morphism $\gamma:B\to L[[t]]$ and an
integer $m>0$, such that 
$$\xymatrix{B\rto^Q\dto_\gamma & K\dto\\ L[[t]]\rto^{t=0} &L}$$
commutes and $\gamma(f)=t^m$. 
\end{lem}
\begin{proof}
Without loss of generality, $B$ is local with maximal ideal $\ker
Q$. Then we can find a discrete valuation ring $A$ inside the quotient
field of $B$ which dominates $B$.  Pass to its completion
$\hat{A}$. The image of $f$ in $\hat{A}$ is of the form $ut^m$, for a
unique $m>0$ and unit $u$, parameter $t$ for $\hat{A}$.  In a suitable
extension $\tilde{A}$ of $\hat{A}$, we can find an $m$-th root of $u$
and change the parameter such that we have that $f$ maps to
$t^m$ in $\tilde{A}$. Choosing a field of representatives $L'$ for
$\tilde{A}$ we get 
an isomorphism $\tilde{A}\cong L'[[t]]$ and hence a morphism
$\gamma':B\to L'[[t]]$ satisfying the requirements of the lemma with
the residue field of $B$ in place of $K$. Passing to a common
extension $L$ of $K$ and $L'$ over this residue field, we obtain
$\gamma$.
\end{proof}

\begin{lem}\Label{paths}
Let $A$ be an integral noetherian $\cc$-algebra and $I\leq A$ an
ideal. Let $Q:\bigoplus_{i\geq0} I^i/I^{i+1}\to K$ be a morphism to a
field, which does not vanish identically on the augmentation
ideal. Then there exists a field extension $L/K$, a morphism
$\gamma:A\to L[[t]]$ and an integer $m>0$ such that the diagram
$$\xymatrix{A\rto\dto_\gamma & A/I \rto^{Q} & K\dto\\
L[[t]]\rrto && L}$$
commutes and 
$$Q(f^{(1)})=\frac{\gamma(f)}{t^m}\Bigr\rvert_{t=0}\,,$$
for every $f\in I$.  Here $f^{(1)}$ denotes the element $f\in I$
considered as an element of the first graded piece of
$\bigoplus_{i\geq0}I^i$.
\end{lem}
\begin{pf}
Choose $g\in I$ such that $Q$ does not vanish on $g^{(1)}$. Apply
Lemma~\ref{value} to
the localization of $\bigoplus_{i\geq0}I^i$ at the element $g^{(1)}$,
the non-zero element $g^{(0)}/g^{(1)}$ and the induced ring morphism
to $K$. We obtain a commutative diagram
$$\xymatrix{\bigoplus_{i\geq0} I^i\rto^-Q\dto_{\tilde{\gamma}} & K\dto\\
L[[t]]\rto &L}$$
and an integer $m>0$
with the property that
$\tilde\gamma(g^{(0)})=t^m\tilde\gamma(g^{(1)})$. We obtain
$\gamma:A\to L[[t]]$ by restricting $\tilde\gamma$ to the degree zero
part of $\bigoplus_{i\geq0}I^i$. Now, for any element $f\in I$ we have
the
equation $f^{(0)}g^{(1)}=g^{(0)}f^{(1)}$ inside
$\bigoplus_{i\geq0}I^i$. Applying $\tilde\gamma$ and cancelling out the unit
$\tilde\gamma(g^{(1)})$, we obtain $\gamma(f)=t^m\tilde\gamma(f^{(1)})$. 
\end{pf}

To prove the theorem, we may assume that $M=\spec A$ is affine and
admits global coordinates $x_1,\ldots, x_n$ giving rise to an \'etale
morphism $M\to\A^n$.  Then we write $\omega=\sum_{i=1}^n f_i\,dx_i$,
for regular functions $f_i$ on $M$.

\begin{lem}\Label{blow}
The conic subscheme $C\subset\Omega_M$ defined by the 1-form
$\omega=\sum_i f_i dx_i$ on $M$ is Lagrangian if 
for every field $K/\cc$, every path $\gamma:\spec K[[t]]\to M$ and
every $m>0$ such that $t^m\mid f_i\big(\gamma(t)\big)$, for all $i$,
we have
\begin{equation}\Label{form}
\sum_{i=1}^n
  d\gamma_i(0)\wedge d\biggl( 
  \frac{f_i\big(\gamma(t)\big)}{t^m}\Bigr\rvert_{t=0} \biggr)=0\,,
\end{equation}
in $\Omega^2_{K/\cc}\,$. Here $\gamma_i=x_i\circ\gamma$. 
\end{lem}
\begin{pf}
First note that as a normal cone, $C$ is pure-dimensional, of
dimension equal to $\dim M$. So 
To prove that $C$ is Lagrangian, we may show that the 2-form
$d\alpha$ defining the symplectic structure on $\Omega_M$ vanishes
when pulled back via $Q:\spec K\to C$, for every morphism $Q$ from the
spectrum of a field to $C$. Moreover, let us note that $d\alpha$ will
vanish on $\spec K$ if it vanishes on $\spec L$ for some extension
$L/K$. 

Note that we have a cartesian diagram
$$\xymatrix{
\Omega_M\rto\dto & \Omega_{\A^n}\dto\\
M\rto & \A^n}$$
The coordinates $p_1,\ldots,p_n$ on $\Omega_{\A^n}=\A^{2n}$ pull back
to functions on $\Omega_M$, which we denote by the same symbols. Thus
$x_1,\ldots,x_n,p_1,\ldots,p_n$ are \'etale coordinates on
$\Omega_M$.  In fact $\Omega_M=\spec A[p_1,\ldots,p_n]$. 
The 2-form
$d\alpha$ is equal to $\sum_i 
dp_i\wedge dx_i$ in these coordinates. 

The ideal defining $X$ is $I=(f_1,\ldots,f_n)\leq A$. The normal cone
$C$ is the spectrum of the graded ring $\bigoplus_{i\geq0}I^i/I^{i+1}$
and the embedding
$C\to \Omega_M$ is given by the ring epimorphism 
$A[p_1\ldots,p_n]\to \sum_{i\geq0}I^i/I^{i+1}$ sending $p_i$ to
$f_i^{(1)}$. 

Thus we have
$$Q^\ast(d\alpha)=Q^\ast\sum dp_i\wedge dx_i=\sum
dQ^\ast(f_i^{(1)})\wedge d Q^\ast(x_i^{(0)})\,.$$
If $Q^\ast$ vanishes on the entire augmentation ideal, this expression
is obviously zero. So assume that $Q$ does not
vanish on the entire augmentation ideal, and choose $\gamma$, $m$ as
in Lemma~\ref{paths}. Then we get
$$Q^\ast(d\alpha)=\sum_{i=1}^n
d\biggl( 
  \frac{f_i\big(\gamma(t)\big)}{t^m}\Bigr\rvert_{t=0} \biggr)
\wedge 
  d\gamma_i(0)\,,
$$
which vanishes by hypothesis. 
\end{pf}

We will now prove the theorem by verifying the condition given in
Lemma~\ref{blow}. Thus we choose 
a field extension $K/\cc$, a path $\gamma:\spec K[[t]]\to M$ and an integer
$m>0$ such that $t^m\mid f_i\big(\gamma(t)\big)$, for all $i$. We
claim that Formula~(\ref{form}) is satisfied in the $K$-vector space
$\Omega^2_{K/\cc}\,$. 

We will introduce some notation.  Define the field elements
$c_i^{(p)}, F^{(p)}_i\in K$ by the formulas
$$\gamma_i(t)=\sum_{p=0}^\infty \frac{1}{p!}c^{(p)}_it^p\,,\qquad
(f_i\circ\gamma)(t)=\sum_{p=0}^\infty\frac{1}{p!}F^{(p)}_i t^p\,.$$
We claim that
\begin{equation}\Label{claim}
\sum_i F_i^{(m)} dc_i^{(0)}=0\,.
\end{equation}
This will finish the proof, because
$$\sum_{i} d\gamma_i(0)\wedge d\biggl(
  \frac{f_i\big(\gamma(t)\big)}{t^m}\Bigr\rvert_{t=0} \biggr)=
-\frac{1}{m!}\, d\sum_i F_i^{(m)} dc_i^{(0)}\,.$$
For future reference, let us remark that the assumption $t^m\mid
  f_i\big(\gamma(t)\big)$, for all $i$, is equivalent to
\begin{equation}\Label{future}
\forall p < m \colon \quad F^{(p)}_i=0\,,
\end{equation}
for all $i$. 

Let us now use the fact that $\omega$ is almost closed.  This means
that 
$$(d\omega)|_X=0\in\Gamma(X,\Omega_M^2|_X)\,.$$ 
By considering the
commutative diagram of schemes
$$\xymatrix{
\spec K[t]/t^m\rto\dto & X\dto \\
\spec K[[t]]\rto^-\gamma & M}$$
we see that this implies that 
\begin{align}\nonumber\Label{property10}
\gamma^\ast(d\omega)|_{\spec
  K[t]/t^m}=0&\in\Gamma(\spec K[t]/t^m,\Omega^2_{K[[t]]}|_{\spec
  K[t]/t^m})\\ &= \Omega^2_{K[[t]]}\otimes_{K[t]}K[t]/t^m\,.
\end{align}
Let us calculate $\gamma^\ast(d\omega)$.  This calculation takes place
inside $\Omega^2_{K[[t]]}\,$: 
\begin{align*}
\gamma^\ast(d\omega)& =\sum_i d(f_i\circ\gamma)\wedge d\gamma_i\\
&=\sum_i\sum_{p=0}^\infty\frac{1}{p!}\left( (dF_i^{(p)})t^p+
F_i^{(p)}pt^{p-1}dt\right)\sum_{p=0}^\infty\frac{1}{p!}\left( (dc_i^{(p)})t^p+
c_i^{(p)}pt^{p-1}dt\right)\\
&= \sum_{p=0}^\infty\frac{1}{p!}\left(
\sum_{k=0}^p {\binom{p}{k}}\sum_i dF_i^{(k)}\wedge
dc_i^{(p-k)}\right)t^p\\
&+ \sum_{p=0}^\infty\frac{1}{p!}\left(
\sum_{k=0}^p {\binom{p}{k}}\sum_i
c_i^{(p+1-k)}dF_i^{(k)}\right)\wedge t^p dt\\ 
&- \sum_{p=0}^\infty\frac{1}{p!}\left(
\sum_{k=0}^p {\binom{p}{k}}\sum_i
F_i^{(k+1)}dc_i^{(p-k)}\right)\wedge t^p dt 
\end{align*}
By Property~(\ref{property10}), the coefficient of $t^{m-1} dt$
vanishes. Thus, the equation
\begin{equation}\Label{holds}
\sum_{k=0}^{m-1}{\binom{m-1}{k}}\sum_ic_i^{(m-k)}dF_i^{(k)}=
\sum_{k=0}^{m-1}{\binom{m-1}{k}}\sum_iF_i^{(k+1)}dc_i^{(m-k-1)}
\end{equation}
holds inside $\Omega_{K/\cc}$. Now, by Assumption~(\ref{future}), all
terms on the left hand side of (\ref{holds}) vanish, as well as the
terms labelled $k=0,\ldots,m-2$ of the right hand side. Hence, the
remaining term on the right hand side of (\ref{holds}) also vanishes.
This is the term labelled $k=m-1$ and is equal to the term claimed to
vanish in~(\ref{claim}).  This concludes the proof of Theorem~\ref{main}.

\subsection{Conclusions}

\subsubsection{Obstruction cones are Lagrangian}

Let $X$ be a Deligne-Mumford stack with a symmetric obstruction
theory.  Suppose $X\hookrightarrow M$ is a closed immersion into a
smooth Deligne-Mumford stack $M$ and let $C\subset \Omega_M$ be the
associated obstruction cone  (see Remark~\ref{roc}).

The following fact was suggested to hold by R. Thomas at the workshop on
Donaldson-Thomas invariants at the University of Illinois at
Urbana-Champaign:

\begin{thm}
The obstruction cone $C$ is Lagrangian.
\end{thm}
\begin{pf}
This follows by combining Theorem~\ref{main} with Corollary~\ref{exists}.
\end{pf}

\begin{cor}\Label{clc}
For the fundamental cycle of the obstruction cone we have
$$[C]=L(\c_X)=\Ch(\nu_X)\,.$$
\end{cor}
\begin{pf}
Because $[C]$ is Lagrangian, we have $[C]=L(\pi[C])$, with
notation as in~(\ref{LL}) and~(\ref{PI}). It remains to show that
$\pi[C]=\c_X$. But this is a local problem, and so we may assume
that our symmetric obstruction theory comes from an almost closed
1-form on $M$. Then $C=C_{X/M}$.
\end{pf}

\subsubsection{Application to Donaldson-Thomas type invariants}

Let $X$ be a quasi-projective Deligne-Mumford stack with a symmetric
obstruction theory.  Let $[X]^\vir$ be the associated virtual
fundamental class. 

\begin{prop}\label{vclass}
We have
$$[X]^\vir=(\alpha_X)_0=c_0^{SM}(\nu_X)\,,$$ 
where $(\alpha_X)_0$ is the degree zero part of the Aluffi class.
\end{prop}
\begin{pf}
Embed $X$ into a smooth Deligne-Mumford stack $M$. Then combine
Proposition~\ref{vireq} with Corollaries~\ref{loc} and~\ref{clc} to
get $[X]^\vir=c_0^M(\c_X)$. 
\end{pf}

\begin{rmk}
In the case that $X=Z(df)$, for a regular function $f$ on a smooth
scheme $M$, the virtual fundamental class is the top Chern class of
$\Omega_M$, localized to $X$. Proposition~\ref{vclass} 
in this case is implicit in
\cite{Alu}. Aluffi proves that $\alpha_X=c(\Omega_M)\cap s(X,M)\in
A_\ast(X)$. Thus, $(\alpha_X)_0=c_n(\Omega_M)\cap[M]\in A_\ast(X)$, by
Proposition~6.1.(a) of \cite{fulton}. 
\end{rmk}

\begin{thm}\Label{finally}
If $X$ is proper, the virtual count is equal to the weighted Euler
characteristic 
$$\#^\vir(X)=\tilde\chi(X)=\chi(X,\nu_X)\,,$$
at least if $X$ is smooth, a global finite group quotient or a gerbe
over a scheme.
\end{thm}
\begin{pf}
Combine Propositions~\ref{vclass} and~\ref{macaluprop} with one
another.
\end{pf}

\begin{rmk}
It should be interesting to prove Theorem~\ref{finally} for arbitrary
proper Deligne-Mumford stacks with a symmetric obstruction theory.
\end{rmk}

\begin{rmk}
Theorem~\ref{finally} applies to all Examples discussed in
Section~\ref{exam} which give rise to proper $X$. In the case of
non-proper $X$, {\em define }$\tilde\chi(X)$ to be the virtual
count. 
\end{rmk}

\begin{rmk}
Let us point out that for a Calabi-Yau threefold the Donaldson-Thomas
and the Gromov-Witten moduli spaces share a large open part, namely
the locus of smooth embedded curves. Both obstruction theories are
symmetric on this locus, and the associated virtual count of this open
locus is the same, for both theories.  

This  observation may or may not be
significant for the conjectures of~\cite{MNOP}.
\end{rmk}

\subsubsection{Another formula for $\nu_X(P)$}

Let $\omega$ be an almost closed 1-form on a smooth scheme $M$. Let
$X=Z(\omega)$ be the scheme-theoretic zero locus of $\omega$ and $P\in
X$ a closed point. Let $x_1,\ldots,x_n$ be \'etale coordinates for $M$
around $P$ and $x_1,\ldots,x_n,p_1\ldots,p_n$ the associated canonical
\'etale coordinates for $\Omega_M$ around $P$.  Write
$\omega=\sum_{i=1}^n f_i dx_i$ in these coordinates. 

Let $\eta\in\cc$ be a non-zero complex number and consider 
the image of the morphism $M\to\Omega_M$ given by the section
$\frac{1}{\eta}\omega\in\Gamma(M,\Omega)$. We call this image
$\Gamma_\eta$. It is a smooth submanifold of $\Omega_M$ of real
dimension $2n$. it is defined by the equations $\eta p_i=f_i(x)$. 

Let $\Delta$ be the image of the morphism $M\to \Omega_M$ given by the
section $d\rho$ of $\Omega_M$, where $\rho=\sum_i x_i\ol{x}_i+\sum_i
p_i\ol{p}_i$ is the square of the distance function defined on
$\Omega_M$ by the choice of coordinates.  Thus $\Delta$ is another
smooth submanifold of $\Omega_M$ of real dimension $2n$. It is defined
by the equations $p_i=\ol{x}_i$. 

Orient $\Gamma_\eta$ and $\Delta$ such that the maps
$M\to\Gamma_\eta$ and $M\to\Delta$ are orientation preserving. 

\begin{prop}
For $\epsilon>0$ sufficiently small, and $|\eta|$ sufficiently small
with respect to $\epsilon$, we have
\begin{equation}\Label{spirit}
\nu_X(P)=L_{S_{\epsilon}}(\Gamma_\eta\cap S_{\epsilon},\Delta\cap
S_{\epsilon}) \,,
\end{equation}
where $S_\epsilon=\{\rho=\epsilon^2\}$ is the sphere of radius
$\epsilon$ in $\Omega_M$ centred at $P$ and $\Gamma_\eta\cap
S_\epsilon$, $\Delta\cap S_\epsilon$ are smooth compact oriented
submanifolds of $S_\epsilon$  with linking number $L$.
\end{prop}
\begin{pf}
Let $C\hookrightarrow\Omega_M$ be the embedding of the normal cone
$C_{X/M}$ into $\Omega_M$ given by $\omega$. Then $\Ch(\nu_X)=[C]$. 
The inverse of $\Ch$ is calculated in Theorem 9.7.11
of~\cite{kash-schap} (see also \cite{Gin}). We get
$$\nu_X(P)=I_{\{P\}}([C],[\Delta])\,,$$
the intersection number at $P$ of the cycles $[C]$ and
$[\Delta]$. Note that $P$ is an isolated point of the intersection
$C\cap\Delta$, by Lemma 11.2.1 of~\cite{Gin}. We should remark that
\cite{kash-schap} deals with the real case.  This introduces various
sign changes, which all cancel out.

Now use Example 19.2.4
in \cite{fulton}, which relates intersection numbers to linking
numbers. We get
$$\nu_X(P)=L_{S_{\epsilon}}([C]\cap S_{\epsilon},\Delta\cap
S_{\epsilon})\,,$$
for sufficiently small $\epsilon$. Next, use Example 18.1.6(d) in
[ibid.], which shows that $\lim_{\eta\to0}[\Gamma_\eta]=[C]$, i.e.,
that there exists an algebraic cycle in $\Omega_M\times \A^1$ which
specializes to $[\Gamma_\eta]$ for $\eta\not=0$ and to $[C]$ for
$\eta=0$. It follows that for sufficiently small $\eta$, we can
replace $[C]$ in our formula by $\Gamma_\eta$.
\end{pf}

\begin{rmk}
Note how Formula~(\ref{spirit}) is similar in spirit to
Formula~(\ref{fp}).  Combining these two formulas for $\nu_X(P)$,
using $\omega=df$, gives
an expression for the  Euler characteristic of the Milnor fibre in
terms of a linking number.
\end{rmk}

\subsubsection{Motivic invariants}

Let $A$ be a commutative ring and $\mu$ an $A$-valued motivic measure
on the category of finite type schemes over $\cc$. For a
scheme $X$, it is tempting to define
$$\tilde\mu(X)=\mu(X,\nu_X)=\int_X\nu_Xd\mu\,$$
and call it the {\em virtual motive }of $X$.  Note that
$\tilde\mu(X)$  encodes the scheme structure of $X$ in a much more
subtle way than the usual motive $\mu(X)$, which neglects all
nilpotents in the structure sheaf of $X$. 

If $X$ is endowed with a symmetric obstruction theory, $\tilde\mu(X)$
may be thought of as a motivic generalization of the virtual count, or
a motivic Donaldson-Thomas type invariant. 

Here are two caveats:

\begin{rmk}
The proper motivic Donaldson-Thomas type invariant should probably
motivate not only $X$ but also $\nu_X$.  For example, in the case of
the singular locus of a a hypersurface, motivic 
vanishing cycles (and not just their Euler characteristics) should
play a role. 
\end{rmk}

\begin{rmk}
Note that $\mu$ will not satisfy Property~(v) of Section~\ref{wec},
unless $\mu=\chi$. So one encounters difficulties when extending the
virtual motive to Deligne-Mumford stacks. To extend $\mu$ to stacks
one formally inverts $\mu(GL_n)$, for all $n$, but then one looses the
specialization to $\chi$, as $\chi(GL_n)=0$. So one cannot think of
the virtual motive of a stack as a generalization of the virtual
count, even if the stack admits a symmetric obstruction theory.
For example, $\tilde\mu(B\zz/2)=1$.
\end{rmk}

\eject


\end{document}